\documentclass[preprint]{imsart}

\pdfoutput=1
\RequirePackage[OT1]{fontenc}
\RequirePackage{amsthm,amsmath,amssymb}
\RequirePackage{natbib}
\usepackage{graphicx}
\usepackage{multirow}
\RequirePackage[colorlinks,citecolor=blue,urlcolor=blue]{hyperref}
\usepackage[margin=1in]{geometry}

\def\T{{ \mathrm{\scriptscriptstyle T} }}

\newtheorem{theorem}{Theorem}

\newtheorem{condition}{Condition}

\usepackage{array,epsfig,fancyheadings,rotating}

\begin{document}

\begin{frontmatter}
\title{Detecting Serial Dependence in Binomial Time Series \textrm{II}: Observation Driven Models}
\runtitle{Detecting Serial Dependence in Binomial Time Series}
\author{W. T. M. Dunsmuir \\School of Mathematics and Statistics, University of New South Wales, Sydney, Australia.\\ w.dunsmuir@unsw.edu.au \\
\and \\
J.Y. He \\School of Mathematics and Statistics, University of New South Wales, Sydney, Australia.\\ jieyi.he@unsw.edu.au}
\maketitle

\begin{abstract}
The detection of serial dependence in binary or binomial valued time series is
difficult using standard time series methods, particularly when there are
regression effects to be modelled. In this paper we derive score-type tests
for detecting departures from independence in the directions of the
GLARMA\ and BARMA\ type observation driven models. These score tests can
easily be applied using a standard logistic regression and so may have appeal
to practitioners who wish to initially assess the need to incorporate serial
dependence effects. To deal with the nuisance parameters in some GLARMA models
a supremum type test is implemented. 
\end{abstract}

\begin{keyword}
\kwd{Binomial time series}
\kwd{Score test}
\kwd{Supremum test}
\kwd{Observation driven models}
\end{keyword}

\end{frontmatter}

\section{Introduction}\label{Sec: Introduction}%

The context in which this paper is relevant is when primary interest is in the detection and estimation of serial dependence in regression models for binomial time series. The need for such a development is clearly demonstrated in an increasing array of applications such as modeling of economic recessions, disease counts, criminal records and sporting events, which are, in most cases, binary or binomial responses. If serial dependence is not detected then standard generalized linear model (GLM) fitting methods can be used to provide correct point estimates  of regression effects and of their standard errors. If serial dependence is detected, depending on the method used, the testing results could provide guidance as to the features of the dependence as a precursor for specifying the form of serial dependence model that might be appropriate.

For the purpose of model development it is useful to have methods that detect
serial dependence without fitting complicated models.
The score test relies only on fitting the model under the assumption that there is no serial dependence and hence is a simple technique for assessing the specification of correlation.
Specifically to assess serial dependence in discrete valued time series, score, or Lagrange multiplier, tests have been developed and applied in previous literature such as \cite{breusch1980lagrange}, \cite{lee1993locally}, \cite{jung2003testing}, \cite{nyberg2008testing} and
\cite{nyberg2010studies}. We propose here the use of score type tests
for testing the null hypothesis that serial dependence is not present against the
alternative that it is induced by an observation driven process.

Let $Y_t$ be a time series taking values in the non-negative integers, $x_{t}$ be an observed $r$-dimensional vector of regressors available at time $t$, $Z_{t}$ a random process and
\begin{equation}\label{eq: LinPredWt}%
W_{t}=x_{t}^{\T}\beta + Z_{t}
\end{equation}
the state variable. Then given the $W_t$, $Y_{t}$ are assumed to be independent with exponential family density which we write in the form
\begin{equation} \label{eq: EFDensity}%
f(y_{t}|W_t)=\exp\left\{y_{t}W_{t}- m_{t}b(W_{t}) + c(y_{t})\right\}.
\end{equation}
Although our methods extend easily and in obvious ways to this general frame work we focus here on the binomial case where $Y_t$ is the number of successes in $m_t$ binomial trials conducted at time $t$, $b(W_t)=\log(1+\exp(W_t))$, $c(y_t)=\log \{m_t!/[y_t!(m_{t}-y_{t})!]\}$, $\mu_t=E(Y_{t}|W_{t})= m_{t}\dot{b}(W_{t})$ and $\sigma_t^2= \mathrm{Var}(Y_{t}|W_{t})= m_{t}\ddot{b}(W_{t})$. Let the probability of a success at time $t$, given $W_{t}$, be denoted
\begin{equation*}
\pi_t =\dot{b}(W_t) = \frac{e^{W_t}}{1+e^{W_t}}, \quad
\sigma_{t}^{2}= m_t\ddot{b}(W_t)=m_{t}\pi_{t}(1-\pi_{t})
\end{equation*}
where $\dot b$ and $\ddot b$ denote first and second derivatives with respect to the argument of $b$.

There are two main specifications of the random process: observation driven where $Z_t$ is specified in terms of previous observations, and parameter driven where $Z_t$ is an unobserved random process. \cite{davis1999modeling} provide and earlier review of these two model structures -- see also \cite{dunsmuir2016} and other articles in the same volume. In this paper we concentrate on methods for testing the null hypothesis that $Z_t$ is absent (no serial dependence)
versus $Z_t$ is an observation driven process. Let $\mathcal{X}_{t}=\{x_{s}: s\le t\}$ and $\mathcal{Y}_{t}=\{y_{s}: s<t\}$, then $Z_t = h(\mathcal{X}_{t},\mathcal{Y}_{t}; \delta)$, where $\delta^{T} = (\beta^T, \psi^T, \omega^T)$ and $(\psi,\omega)$ are the parameters specifying the particular form of $Z_t$. Here $\psi$ are the parameters which, if set to zero, give $Z_t \equiv 0$ (subject to suitable initial conditions) and under this hypothesis $\omega$ are nuisance parameters that are not estimable.

Observation driven models take various forms - - see \cite{BenjaminJASA2003} for a
general discussion. We consider two classes of observation driven models
here: generalized linear autoregressive moving average (GLARMA)\ models reviewed in general in  \cite{dunsmuir2016} and for the Poisson case in \cite{davis2000ObsDriven}, \cite{davis2005maximum} and binary autoregressive moving average (BARMA)\ models given in \cite{wang2011autopersistence}.
In this paper we will derive the score statistics tested against
the GLARMA and BARMA models with binomial responses under the null hypothesis.
In a companion paper, \cite{dunsmuir2016testing}, we derive the
score-type test for detecting parameter driven serial dependence where $Z_t$ is a latent stationary random process. These are of substantially different structural form and require different large sample theory to that considered here.

The outline of this paper is as follow. In Section \ref{Sec: GLARMA models} we review the GLARMA\ model.
As part of this review we derive the relevant likelihoods, score functions and
information matrices for these models and develop the asymptotic theory for
score statistics for testing the null hypothesis of no lagged dependence terms
in the state equations. We also derive the supremum score test statistic
for GLARMA models with nuisance parameters and investigate its asymptotic
distribution for simple examples. Section \ref{Sec: BARMA models} reviews the BARMA\ model and its score functions. Section \ref{Sec: Other tests} introduces alternative tests for serial dependence such as the Box-Pierce-Ljung, likelihood ratio and Wald tests. Section \ref{Sec: simulation} assesses the asymptotic results of supremum score test using finite sample simulations. Section \ref{Sec: applications} applies these ideas to some real data series. Section \ref{Sec: conclusions} concludes.

\section{GLARMA\ Models} \label{Sec: GLARMA models}%

Written in a general format, the `noise' process of a GLARMA model is
\begin{equation}\label{eq: GLARMA Zt}%
Z_{t}=\sum_{j\in J_{\phi}}\phi_{j}\left( Z_{t-j}+e_{t-j}\right)+ \sum_{j\in J_{\theta}}
\theta_j e_{t-j}
\end{equation}%
where $J_{\phi}$ is the set of autoregressive lags with non-zero $\phi_j$ and $J_{\theta}$ is the set of moving average lags with non-zero $\theta_j$.

In \eqref{eq: GLARMA Zt} the ``residuals", $e_t$, can be defined in a number of ways. Starting with what we refer to as identity residuals, $e_{t}^{I}=y_{t}-m_t \pi_{t}$, we consider
\begin{equation}\label{eq: GammaResids}
e_{t}=\sigma_{t}^{-\gamma}e_{t}^{I}
\end{equation}
where $\gamma=0$ gives identity residuals introduced by \cite{wang2011autopersistence}, $\gamma=1$ Pearson residuals and $\gamma=2$ ``score-type" residuals used in \cite{creal2008general} -- these are the three types of residuals commonly encountered and are supported in the R-language \texttt{glarma} package described in \cite{dunsmuir2014glarma}. Note that $E(e_t)=0$ but only the Pearson residuals have unit variance.

The recursion for $Z_t$ in \eqref{eq: GLARMA Zt} can be rewritten as
\begin{equation}\label{eq: GLARMA ZtAlt}%
Z_{t}=\sum_{j\in J_{\phi}\bigcap J_{\theta}}\omega_{j} Z_{t-j} + \sum_{j\in J_{\phi}\bigcap J_{\theta}^C}\psi_{j} Z_{t-j} + \sum_{j \in J_{\phi}\bigcup J_{\theta}}\psi_j e_{t-j}
\end{equation}
where $\omega_{j} = \phi_{j}$ and $\psi_j=\theta_j+\phi_j$ for $j\in J_{\phi}\bigcap J_\theta$; $\psi_j = \phi_j$ for $j\in J_{\phi}\bigcap J_{\theta}^C$ and $\psi_j=\theta_j$ for $j\in J_{\phi}^C\bigcap J_\theta$. Since $J_\phi \bigcap J_\theta^C \subset J_\phi\bigcup J_\theta$, the null hypothesis $H_0:\psi=0$ specifies that the second and third summations in \eqref{eq: GLARMA ZtAlt} are zero and \eqref{eq: GLARMA ZtAlt} reduces to a recursive equation in
$Z_t = \sum_{j\in J_\phi\bigcap J_\theta} \omega_j Z_{t-j}$ with $0$ as its solution (provided the initial conditions are zero). The $\omega$ are nuisance parameters that cannot be estimated under the null hypothesis. The nuisance parameters can be taken into account using the general tests proposed by \cite{davies1977hypothesis}, \cite{davies1987hypothesis} and \cite{andrews1994optimal}.

\subsection{Score tests against GLARMA alternatives}\label{SSc: ST GLARMA}%

Throughout we assume that the recursion \eqref{eq: GLARMA ZtAlt} is initialized using pre-observation period values $Z_t = e_t = 0$ for $t \le 0$, that is, setting them to their unconditional stationary mean values. Under this assumption the (conditional) log-likelihood is
\begin{equation}\label{eq: GLARMA log likelihood}
l(\delta)= \sum_{t=1}^n Y_tW_t(\delta)-m_t b(W_t(\delta))+ c(y_t).
\end{equation}
Let $\delta_{0}=(\beta_{0}, 0, \omega)$ be the true parameter value under $H_{0}:\psi=0$, since $\omega$ are the nuisance parameters that are not estimable under the null, throughout the derivations of score vectors we assume $\omega$ to be fixed. For fixed $\omega$, denote $\hat\delta_{0}=(\hat \beta_{0}, 0, \omega)$ as the m.l.e.of \eqref{eq: GLARMA log likelihood} under the null, where $\hat \beta_{0}$ is the GLM estimate. Similarly, $\hat\delta = (\hat\beta, \hat\psi, \omega)$ is m.l.e.of
\eqref{eq: GLARMA log likelihood} under the alternative, and $\omega$ is fixed.

The score vector is $S(\delta)=\partial l(\delta)/\partial \delta$ which requires recursive calculation with $W_t(\delta)$, $Z_t(\delta)$, $e_t(\delta)$, $\sigma_t(\delta)$ and their derivatives. We assume these recursions, as well as their derivatives (e.g. $\partial Z_t/\partial \delta$,
$\partial e_t/\partial \delta$), are also initialized at zero for $t \le 0$.

Denote $J_L=J_\phi\bigcup J_\theta =\{j_1,\ldots, j_{L}\}$ and $e_{t-J_L}(\delta) =(e_{t-j_1}(\delta), \ldots, e_{t-j_L}(\delta))^\T$. Now
\begin{equation}\label{eq: Deriv1st GLARMA}%
\frac{\partial Z_{t}}{\partial\psi}\vert_{\delta_0} = e_{t- J_L}(\delta_0) +
\sum_{j\in J_\phi\bigcap J_\theta} \omega_{j}\frac{\partial Z_{t-j}(\delta_0)}{\partial\psi}.
\end{equation}
Rewrite $\left(1-\sum_{j\in J_\phi\bigcap J_\theta} \omega_{j}\xi^{j}\right)^{-1}e_{t-a} =
\sum_{j=0} ^{\infty}\tau_{j}(\omega)e_{t-a-j}$ for $a\in J_L$, where $\tau_{0}(\omega) =1$,
so that solving \eqref{eq: Deriv1st GLARMA},
\begin{equation}
\frac{\partial Z_{t}}{\partial\psi}\vert_{\delta_0}=\sum_{j=0} ^{\infty}\tau_{j}(\omega)e_{t-J_L-j}(\delta_{0}).
\end{equation}
Also, it is straightforward to show that $\partial Z_{t}/\partial \beta|_{\delta_0} = 0$,
$\partial Z_{t}/\partial \omega|_{\delta_0}=0$. Thus for any fixed $\omega$, the scaled
score vector is
\begin{equation}\label{eq: ST GLARMA Gen}%
S(\delta_{0}) = n^{-1/2} \sum_{t=1}^{n}\left( y_{t}-m_t \pi_{t}(\delta_{0})\right)
\begin{bmatrix}  x_{t}\\
\sum_{j=0} ^{\infty}\tau_{j}(\omega)e_{t-J_L-j}(\delta_0)\\
0 \end{bmatrix}.
\end{equation}
The covariance matrix corresponding to the component of the score vector for $\beta$ and $\psi$ is
\begin{equation}\label{eq: InfMat GLARMA Gen}%
I_{n}(\delta_0) = n^{-1}\sum_{t=1}^n \sigma_t^2(\beta_0)
\begin{bmatrix}
x_{t}x_{t}^{\T} & 0 \\
0 & \Gamma_{t,L}(\delta_0)
\end{bmatrix}
\end{equation}
and $\Gamma_{t,L}$ is $L\times L$ symmetric matrix
\begin{equation*}
\Gamma_{t,L} = \sum_{h=0}^{(L-1)}\sum_{j=0}^{\infty}\Gamma_{t-j, h},
\end{equation*}
\begin{equation*}
\Gamma_{t-j,0} = \tau_{j}(\omega)^2\textrm{diag}\left[(\sigma_{(t-j_{1})-j}^{2-2\gamma}, \ldots, \sigma_{(t-j_L)-j}^{2-2\gamma})\right],
\end{equation*}

\begin{eqnarray*}
\Gamma_{t-j, h} = A_{t-j, h} + A^{T}_{t-j, h}, && h=1,\ldots,L-1
\end{eqnarray*}
and
\begin{equation*}
A_{t-j,h}=\textrm{h-superdiag} \left[\tau_j(\omega)\tau_{j+(j_{h+1}-j_1)}(\omega) \sigma_{(t-j_{h+1})-j}^{2-2\gamma}, \ldots, \tau_{j}(\omega)\tau_{j+(j_L-j_{L-h})}(\omega) \sigma_{(t-j_L)-j}^{2-2\gamma}\right].
\end{equation*}

In practice, the score vector $S(\delta_0)$ is evaluated by replacing $\delta_0$ with
$\hat\delta_0$. By the definition of $\hat\beta_0$,
$\sum_{t=1}^{n}(y_t-m_t\pi_t(\hat\beta_0))x_t=0$ so that the only non-zero part of $S(\hat \delta_{0})$ in \eqref{eq: ST GLARMA Gen} is
\begin{equation} \label{eq: GLARMA StAlt General}%
S(\hat\delta_0) = n^{-1/2}\sum_{t=1}^{n}\left( y_t -m_t\pi_t(\hat\delta_0)\right)
\left(\sum_{j=0} ^{\infty}\tau_{j}(\omega)e_{t-J_L-j}(\hat\delta_0)\right)
\end{equation}
and its covariance matrix estimated at $\hat \delta_{0}$ is
\begin{equation}\label{eq: InfMat GLARMA Red}%
I_{L}(\hat \delta_0) = n^{-1}\sum_{t=1}^n \sigma_{t}^2(\hat \beta_0)\Gamma_{t,L}(\hat\delta_0).
\end{equation}
The resulting score statistic is
\begin{equation}\label{eq: QStatST}%
Q_{L}^{ST}(\omega) = S(\hat{\delta}_{0})^{\mathrm{T}} I_{L}(\hat\delta_0)^{-1} S(\hat\delta_{0}).
\end{equation}

The consideration of nuisance parameters $\omega$ complicates the test statistics. A simple approach is to fix the value of the nuisance parameters, for example, set $\omega=0$. This results in the same test statistics as the situation when the lags in $\phi$ do not overlap with lags in $\theta$. In this situation $\psi=0$ implies that $\omega=0$ also and there are no nuisance parameters under the null hypothesis of no serial dependence. An example of this situation is when $J_\phi=\{1\}$ and $J_\theta=\{2\}$ then $J_{\phi}\bigcap J_{\theta} =\emptyset$ and
\begin{equation*}
Z_t = \phi (Z_{t-1} + e_{t-1}) + \theta e_{t-2} = \phi Z_{t-1} + \phi e_{t-1} + \theta e_{t-2}.
\end{equation*}
In these situations, the model \eqref{eq: GLARMA ZtAlt} has $\omega$ a subvector of $\psi$. Since under $H_{0}:\psi=0$, $\omega = 0$ also, and the score vector \eqref{eq: QStatST} simplifies to
\begin{equation*}
S(\hat \beta_{0})= n^{-1/2} \sum_{t=1}^{n}(y_t-m_t\pi_t(\hat \beta_0))e_{t-J_L}(\hat \beta_0)
\end{equation*}
with covariance matrix in \eqref{eq: InfMat GLARMA Red} simplifying to%
\begin{equation}\label{eq: InfMat1 GLARMA H0}%
I_{L}(\hat\beta_0) = n^{-1}\sum_{t=1}^n \sigma_{t}^2(\hat\beta_0)\cdot \textrm{diag}\left( \sigma_{t-j_1}^{2-2\gamma}(\hat\beta_0), \ldots, \sigma_{t-j_L}^{2-2\gamma}(\hat\beta_0) \right).
\end{equation}
The resulting score statistic \eqref{eq: QStatST} is
\begin{equation}\label{eq: STGLARMA Std}%
Q_{L}^{ST}(0) = S(\hat\beta_0)^T I_{L}(\hat\beta_0)^{-1}S(\hat\beta_0) = \sum_{l=1}^{L}\hat{C}^2(l)/\hat{B}(l)
\end{equation}
in which
\begin{equation*}
\hat{C}(l) = n^{-1}\sum_{t=j_l+1}^{n}\sigma_{t-j_l}^{-\gamma}(\hat\beta_0)e^{I}_{t}(\hat \beta_0) e^{I}_{t-j_l}(\hat \beta_{0}), \quad \hat{B}(l) = n^{-1}\sum_{t=j_l+1}^n \sigma_t^2(\hat{\beta}_0) \sigma_{t-j_l}^{2(1-\gamma)}(\hat{\beta}_0)
\end{equation*}
where $\gamma = 0, 1, 2$ corresponds to Identity, Pearson and Score residuals \eqref{eq: GammaResids}, respectively, in the GLARMA specification \eqref{eq: GLARMA Zt}. This score statistic is the same for testing against the alternative that the model is a pure AR$(L)$ or a pure MA$(L)$ and hence is a pure
significance test as noted in \cite{poskitt1980testing} for score test of ARMA$(p,q)$ against
ARMA$(p+r, q+s)$.

A better approach, but which is more complicated in its implementation and derivation of asymptotic properties, is the supremum test method proposed in \cite{davies1977hypothesis} and \cite{davies1987hypothesis}. The essential idea of this supremum test statistic is to take the maximum value of the test statistic \eqref{eq: QStatST} over a suitably chosen subset, $\Omega$, of the nuisance parameter space to get
\begin{equation}\label{eq: sup QStatST}
\underset{\omega\in \Omega}\sup~Q_{L}^{ST}(\omega) = \underset{\omega\in \Omega}\sup
S(\hat\delta_0) I_{L}(\hat\delta_0)^{-1}S(\hat\delta_0).
\end{equation}
The asymptotic distribution of supremum tests with a single valued nuisance parameter has been investigated -- see \cite{davies1977hypothesis} and \cite{davies1987hypothesis}, hence to assess
the accuracy of the asymptotic property of supremum version of score test, in simulation
we use examples where the nuisance parameter is one dimensional with space
$\Omega =[\omega_{\mathcal{L}},\omega_{\mathcal{U}}]$. \par

\subsection{Asymptotic properties of GLARMA score statistics}\label{SSc: Asym GLARMA}%

In order to establish the asymptotic distribution of the score test and the likelihood ratio and Wald tests (considered in Section 4.2) the large sample properties of the GLM and GLARMA estimators are required under the null hypothesis. Theorem \ref{Thm: Aysm GLARMA Nui} gives the result for GLARMA estimators which also gives an obvious corollary for the asymptotic properties of the GLM estimators
(separate proof of the asymptotic normality of GLM estimators for binomial responses is provided in \cite{dunsmuir2016marginal}). Some regularity conditions are required:

\begin{condition}\label{cond: mt}
The sequence of trials $\{m_{t}: 1\le m_{t}\le M\}$ is specified in one of two ways:
\begin{description}
\item (a) A stationary process independent of the regressors $\{X_{t}\}$ with $\kappa_j= P(m_t=j)$, $\kappa_M>0$, $\sum_{j=1}^{M}\kappa_j = 1$.

\item (b) A deterministic sequence which are asymptotically stationary and for which $\kappa_{j}$ are limits of finite sample sequences of $m_{t}=j$.
\end{description}
\end{condition}

\begin{condition}\label{cond: Reg Trend Type}
The regression sequence is specified in one of two ways:
\begin{description}
		\item (a) Deterministic covariates defined with functions:  $x_{nt}=h(t/n)$ for some specified piecewise continuous vector function $h: [0,1]\to \mathbb{R}^r$, or,
		
		\item (b) Stochastic covariates which are a stationary vector process: $x_{nt}=x_t$ for all $n$ where $\{x_t\}$ is an observed trajectory of a stationary process for which $E(e^{s^T X_t}) <\infty$ for all $s\in \mathbb{R}^r$.
	\end{description}
\end{condition}

\begin{condition}\label{cond: Xt}
For any fixed $\beta\in \mathbb{R}^r$, as $n\to\infty$, $\underset{1\le t\le n}\sup \Vert n^{-1/2} x_{t} e^{x_t^{\T}\beta}\Vert \overset{p}\to 0$, and the parameter space $\mathbb{X}=\{x_{t}: 1\le t\le n\}$
has $\texttt{rank}(\texttt{span}(\mathbb{X}))=r$.
\end{condition}

The full rank assumption is needed to maintain the non-singularity of the information matrix, thus the consistency and asymptotic normality of observation driven model estimators can be achieved.

\begin{condition}\label{cond: omega}%
$\vert 1- \underset{j\in J_{\phi}\bigcap J_{\theta}}\sum \omega_j\xi^j \vert \ne 0$ for all
$\vert \xi\vert\le 1$. $\xi$ is the backshift operator.
\end{condition}

\begin{theorem}\label{Thm: Aysm GLARMA Nui}
Given Conditions \ref{cond: mt} to \ref{cond: omega}, under $H_{0}:\psi=0$, for fixed $\omega$, as $n\to\infty$, $\hat\delta \to \delta_{0}$ in probability and
$n^{1/2}((\hat \beta, \hat\psi) - (\beta_{0}, 0))\to N(0,I(\delta_{0}))$ in distribution, where
$I(\delta_0) = \underset{n\to\infty}\lim I_{n}(\delta_0)$ and $I_{n}(\delta_0)$ is defined in
\eqref{eq: InfMat GLARMA Gen}.
\end{theorem}
We now use this theorem, applied to the GLM estimators, to get the asymptotic chi-squared distribution for the score statistic.
\begin{theorem} \label{Thm: QNui GLARMA LimDist}%
Assume $\omega$ is fixed, under Conditions \ref{cond: mt} to \ref{cond: omega}, for any fixed $L$ the score statistics $Q_{L}^{ST}(\omega)$ from model \eqref{eq: QStatST} for testing against GLARMA\ alternatives has an asymptotic $\chi^{2}_{L}$ distribution.
\end{theorem}
Theorem \ref{Thm: QNui GLARMA LimDist} covers the case when $\omega=0$ and $Q_{L}^{ST}(\omega) = Q_{L}^{ST}(0)$ in \eqref{eq: STGLARMA Std}.

\subsection{Asymptotic distribution of the supremum score test}\label{SSc: Asym Supremum}%

\cite{davies1987hypothesis} proposed an upper bound for the upper tail probability of the supremum score statistic and $\omega$ is the one dimensional nuisance parameter:
\begin{equation} \label{eq: Sup chi2 AsyDist}%
P\left\{\underset{\omega \in \Omega}\sup S(\omega) > u \right\} \le P(\chi^2_{s} > u) + \int_{\omega_{\mathcal{L}}}^{\omega_{\mathcal{U}}}\psi(\omega)d\omega
\end{equation}
in which $S(\omega)= Z_{1}^2(\omega) + \ldots + Z_{s}^2(\omega)$ where $Z_{i}(\omega)\sim N(0,1)$ for all $i=1,\ldots,s$, and
\begin{equation*}
\psi(\omega) = \frac{1}{\sqrt{2\pi}}\int_{0}^{\infty} \left\{ 1- \prod_{j=1}^{s}(1+\lambda_{j}(\omega)t)^{-1/2} \right\} t^{-3/2}dt \cdot u^{\frac{s-1}{2}}e^{-\frac{u}{2}} \pi^{-\frac{1}{2}}2^{-\frac{s}{2}}/\Gamma(\frac{s}{2}+\frac{1}{2})
\end{equation*}
where $\lambda_{j}(\omega)$, $j=1,\ldots,s$ are the eigenvalues of the matrix $B(\omega)- A^{\T}(\omega)A(\omega)$. Here $Y(\omega) = \partial Z(\omega)/\partial \omega$,
\begin{equation*}
\mathrm{Var}\binom{Z(\omega)}{Y(\omega)} =
\begin{bmatrix}
I & A(\omega)\\
A^{\T}(\omega) & B(\omega)
\end{bmatrix}.
\end{equation*}

We illustrate the application of \cite{davies1987hypothesis} to the simple case where $s=1$,
\begin{equation} \label{eq: Prob upcrossing dim one}%
\int_{\omega_{\mathcal{L}}}^{\omega_{\mathcal{U}}}\psi(\omega)d\omega = \pi^{-1} e^{-\frac{u}{2}} \int_{\omega_{\mathcal{L}}}^{\omega_{\mathcal{U}}}\lambda^{1/2}(\omega)d\omega.
\end{equation}
For the multi-dimensional score vectors of $s\ge 2$, the above integral can be obtained with the
same method but requires the evaluation of higher dimensional integrals.

We next derive the specific details for the GLARMA$(1,1)$ model with Pearson residuals. Suppose
\begin{equation*}\label{eq: GALRMA(1,1) ZtAlt}%
Z_{t}=\phi_{1}(Z_{t-1}+ e_{t-1}^P)+\theta_1 e_{t-1}^P = \omega Z_{t-1} + \psi e_{t-1}^P,
\quad \omega\in (-1, 1).
\end{equation*}
The score vector evaluated under the null is
\begin{equation}\label{eq: GLARMA(1,1) ST}%
S(\delta_0)= n^{-1/2}\sum_{t=1}^n (y_t - m_t\pi_{t}(\beta_0)) \left(\sum_{i=0}^{(t-2)} \omega^{i}e_{t-1-i}^P(\beta_0)\right).
\end{equation}

Let $\gamma^2(\omega)=\mathrm{Var}\left(S(\delta_0)\right)$, based on
Theorem \ref{Thm: QNui GLARMA LimDist} the standardization of the scaled score vector $Z(\omega)=S(\delta_{0})/\gamma(\omega)$ is asymptotically normally distributed with unit
covariance, thus the distribution of $\sup Z^2(\omega)$ follows
\eqref{eq: Sup chi2 AsyDist}, in which $\mathrm{Cov}(Z(\omega), Y(\omega))=0$,
and $\lambda(\omega)=\mathrm{Var}(Y(\omega))\approx (1-\omega^2)^{-2}$ as $n\to\infty$.
By Theorem \ref{Thm: QNui GLARMA LimDist}, $Q_{1}^{ST}(\omega)\to Z^2(\omega)$ in distribution for any fixed $\omega$, so that the distribution in \eqref{eq: Sup chi2 AsyDist} can be rewritten as
$P\left\{\sup_{\omega\in \Omega} Q_{1}^{ST}(\omega) > u\right\}\le \mathcal{F}_{\Omega}(u)$ where
\begin{equation}\label{eq: sup chi-square1 dist}%
\mathcal{F}_{\Omega}(u) = P(\chi^2_1 > u) + \frac{1}{2\pi} e^{-\frac{u}{2}} \ln \left. \left[\frac{1+\omega}{1-\omega} \right] \right\vert_{\omega_{\mathcal{L}}} ^{\omega_{\mathcal{U}}}
\end{equation}
To assess the utility of this upper tail bound we will compare, via simulation, the quartiles of $F_{\Omega}(u)$ with the empirical quantiles of the supremum score test. \par

\section{BARMA\ Models} \label{Sec: BARMA models}%

The recent paper by \cite{wang2011autopersistence} considers the BARMA\ model in which the
model of serial dependence $Z_{t}$ is defined as
\begin{equation}\label{eq: BARMA Zt}
Z_{t}=\sum_{i=1}^{p}\phi_{i}Y_{t-i}+\sum_{i=1}^{q}\theta_{i}e_{t-i}^{I}%
\end{equation}
using the unstandardised residuals. Generalisation of the BARMA\ model to
include scaled residuals such as the Pearson or score residuals
introduced for the GLARMA\ model does not seem to be a sensible idea since
the scale of $Y_t$ would be different from that of the scaled residuals. The BARMA\
model can also be generalized by using $J_{\phi}$ as the set of lags of past observations and $J_\theta$ as the set of lags for residuals. BARMA\ models cannot be written in the GLARMA\
form except when $J_\phi=\emptyset$ and the residuals are specified as $e^{I}_{t-j}$
in both models.

The identifiability issue could also rise for BARMA\ model under the null hypothesis but to a very limited extent. An alternative expression of model \eqref{eq: BARMA Zt} is
\begin{equation*}
Z_{t} = \sum_{j\in J_\phi\bigcup J_\theta}\tilde\phi_{j}Y_{t-j} - \sum_{j\in J_\phi\bigcap J_\theta}\theta_{j}\pi_{t-j} + \sum_{j\in J_\phi^C\bigcap J_\theta}\theta_{j} e_{t-j} ^{I},
\end{equation*}
the last part is zero if $J_\theta \subseteq J_\phi$. In the above function,
$\tilde{\phi}_{j} = \phi_{j}+\theta_{j}$ for $j\in J_\phi\bigcap J_\theta$,
and $\tilde{\phi}_{j} = \phi_{j}$ if $j\in J_\phi\bigcap J_\theta^C$.
If the number of trials are not time varying ($m_t=m$) and the means $\{ \pi_t\}$
are constants (no exogenous covariates), assume $\tilde{\phi} = 0$ and
$\theta_{j} = 0$ for $j\in J_\phi^C\bigcap J_\theta$, then $Z_{t} = 0$ when
$\sum_{j\in J_\phi\bigcap J_\theta}\theta_{j} = 0$ but not all $\theta_{j}$
should be zero. Also, it is easy to show that the BARMA\ model with constant
trials and regressors has non-invertible information matrix under the true parameter.
In the following discussions we exclude this example and focus on time varying regressors only.

Let $\delta=(\beta,\phi,\theta)$ be the parameters in BARMA.
Unlike the GLARMA\ model, the BARMA\ model has less restrictions on its coefficients
$(\phi, \theta)$, as $\{Z_{t}\}$ is bounded uniformly under binomial responses.
Here we are testing that $H_{0}:(\phi, \theta)=0$ versus the alternative
$H_{a}: (\phi,\theta)\ne0$. We will also consider testing that there is no AR
part, no MA part (i.e. separate tests for $\phi=0$ and $\theta=0)$.
Similarly to the GLARMA\ model we have%
\begin{equation}\label{eq: Deriv1st BARMA}%
\frac{\partial Z_{t}}{\partial\phi}= Y_{t-J_\phi} +\sum_{j\in J_\theta}\theta_{j}
\frac{\partial e_{t-j}^I}{\partial\phi}, \quad
\frac{\partial Z_{t}}{\partial\theta}= e_{t-J_\theta} +\sum_{j\in J_\theta}\theta_{j}
\frac{\partial e_{t-j}^I}{\partial\theta}.%
\end{equation}
The derivative at the true value $\delta_0$, under the null, gives a score vector of%
\begin{equation*}
S(\delta_0)=\sum_{t=1}^{n}( y_{t}- m_t\pi_{t}(\delta_0))
\begin{bmatrix}
0\\
Y_{t-J_{\phi}}\\
e_{t-J_{\theta}}^I(\delta_0)
\end{bmatrix}.
\end{equation*}
Using above derivatives of $Z_{t}$, the information matrix for the BARMA\ model is \begin{equation}\label{eq: InfMat BARMA H0}%
I_{n}(\delta_{0}) =\sum_{t=1}^{n} \sigma^2_t(\delta_{0})%
\begin{bmatrix}
x_{t}x_{t}^\mathrm{T} & x_{t}(m\pi)_{t-J_\phi}^\mathrm{T} & 0\\
(m\pi)_{t-J_\phi}x_{t}^\mathrm{T} & A_{\phi\phi,t} & A_{\phi\theta,t}\\
0 & A_{\theta\phi,t} & A_{\theta\theta,t}%
\end{bmatrix}
\end{equation}
in which $(m\pi)_{t-J_{\phi}} = (m_{t-j_1}\pi_{t-j_1}, \ldots, m_{t-j_p}\pi_{t-j_p})$.
$A_{\phi\phi,t}=E(Y_{t-J_\phi}Y_{t-J_\phi}^\mathrm{T}) = \textrm{diag}((m\pi)_{t-J_\phi}(\delta_0))$ and off diagonal elements $m_{t-a}\pi_{t-a}(\delta_0)m_{t-b}\pi_{t-b}(\delta_0)$; $A_{\theta\theta,t}= \textrm{diag}(\sigma^2_{t-J_\theta}(\delta_0))$; thus $A_{\phi\phi,t} = \textrm{diag}(\sigma_{t-J_\phi}^2(\delta_0)) +\Pi_{J_\phi,t}\Pi_{J_\phi,t} ^{\mathrm{T}}$,
$\Pi_{J_\phi,t}= (m\pi)_{t-J_\phi}(\delta_0)$.
\begin{equation*}
A_{\phi\theta,t}=V_{[J_\phi],[J_\theta],t}, \quad A_{\phi\theta,t}= V_{[J_\theta],[J_\phi],t},
\end{equation*}
where $V_{[a],[b],t}$ is a $[a]\times [b]$ matrix. $[a]$ indicates the cardinality of set $a$, $a=\{a_1, \ldots, a_p\}$ and $a_1\le \cdots \le a_p$ (the same for $[b]$). The rows of $V_{[a],[b],t}$ are labeled by components of set $a$; the columns are labeled by components of set $b$. The intersection labeled by $l=a\bigcap b$ have value $\sigma_{t-l}^2(\delta_0)$, otherwise zero.

Denote $\hat{\delta}_0$ as the estimates of $\delta_{0}=(\beta_{0},0,0)$ under the null. To construct the score statistic we partition the information matrix evaluated at $\hat\delta_{0}$ as
\begin{equation*}
I_{n}(\hat{\delta}_0)=
\begin{bmatrix}
\hat{E}_n & \hat{F}_n^\mathrm{T} \\
\hat{F}_n & \hat{G}_n
\end{bmatrix}
\end{equation*}
where
\begin{equation}\label{eq: En}%
\hat{E}_n = \sum_{t=1}^{n}\sigma^2_t(\hat{\delta}_{0}) x_{t}x_{t}^{\mathrm{T}},
\end{equation}
\begin{equation}\label{eq: Fn}%
\hat{F}_n^\mathrm{T} = \sum_{t=1}^{n}\sigma^2_t(\hat{\delta}_{0})
\begin{bmatrix}
x_{t}(m \hat{\pi})_{t-J_\theta}^{\mathrm{T}} & 0
\end{bmatrix},
\end{equation}
\begin{equation}\label{eq: Gn}%
\hat{G}_n^\mathrm{T} = \sum_{t=1}^{n}\sigma^2_t(\hat{\delta}_{0})
\begin{bmatrix}
A_{\hat{\phi}\hat{\phi},t} & A_{\hat{\phi}\hat{\theta},t}\\
A_{\hat{\theta}\hat{\phi},t} & A_{\hat{\theta}\hat{\theta},t}
\end{bmatrix}.
\end{equation}

As before, the component of the score vector corresponding to $\beta$ evaluated at
$\hat{\beta}_0$ is zero so the score statistic for testing for serial dependence is
\begin{equation}\label{eq: QStatBARMA}
Q_L^B(\hat\delta_0) = \begin{bmatrix}
S_{\phi}(\hat \delta_{0}) & S_{\theta}(\hat{\delta}_{0})
\end{bmatrix}
[\hat{G}_n-\hat{F}_n \hat{E}_n^{-1} \hat{F}_n^\mathrm{T} ]^{-1}
\begin{bmatrix}
S_{\phi}(\hat{\delta}_{0}) & S_{\theta}(\hat \delta_{0})
\end{bmatrix}.
\end{equation}
Unlike the score tests against GLARMA\ alternatives this statistic does not simplify to the form of sum of squares of weighted estimates of autocorrelations based on identity residuals.

\begin{theorem}\label{Thm: Asym ST BARMA}%
Under Conditions \ref{cond: mt} to \ref{cond: Xt}, for any fixed $L$ the score statistic $Q_{L}^{B}(\hat \delta_0)$ for testing against BARMA\ alternatives has an asymptotic $\chi^2_{L}$ distribution.
\end{theorem}

\section{Other Test Statistics}\label{Sec: Other tests}

\subsection{Test based on the autocorrelation of Pearson residuals}

The Box-Pearce-Ljung test is based on the Pearson residuals
$e_{t}^{P}(\hat\beta_{0}) = \sigma_t(\hat\beta_0)^{-1}e_{t}^{I}(\hat\beta_0)$,
where $\hat\beta_{0}$ is the GLM estimate and $\hat W_t = x_t^\T \hat\beta_0$.
The auto-covariances are defined, in the usual way, as
\begin{equation} \label{eq: ACVF Pearson}%
C(l)=\frac{1}{n}\sum_{t=1}^{n} e_{t}^{P}(\hat{\beta}_0)e_{t-l}^{P}(\hat{\beta}_0),
\end{equation}
giving estimated autocorrelations as $r(l)=C(l)/C(0)$. The usual Box-Pearce-Ljung
statistic based on these autocorrelations is defined as
\begin{equation} \label{eq: QStatBLP}%
Q_{L}^{BLP}=n(n+2) \sum_{l=1}^{L} (n-l)^{-1} r(l)^2.
\end{equation}
Using Theorem \ref{Thm: QNui GLARMA LimDist}, $Q_{L}^{BLP}$ has an asymptotic $\chi^2_{L}$ distribution in the same way as $Q_{L}^{ST}(0)$. Note that none of the score test statistics based on the three types of residuals correspond to the
$Q_{L}^{BLP}$.

\subsection{\bf 4.2. Tests based on the likelihood estimates}

\noindent {\bf 4.2.1. Likelihood ratio test}

Based on the the log-likelihood given in \eqref{eq: GLARMA log likelihood}, the likelihood ratio statistic for testing the total of $L$ parameters specifying the serial dependence is
\begin{equation}\label{eq: QStatLR}%
Q_{L}^{LR}= 2\left[ l_n(\hat{\delta})-l_n(\hat{\delta}_{0})  \right].
\end{equation}
For regular cases that do not have nuisance parameters (e.g. $Z_t$ is a pure AR or MA process, or there is no overlap between the lags of $\phi$ and $\theta$), the likelihood ratio test statistic
has a large sample chi-square distribution, with degrees of freedom given by the length of correlation parameter $(\phi,\theta)$, under the null hypothesis of no serial dependence. When there are nuisance parameters which can not be estimated under the null, the likelihood ratio test does not have its standard asymptotic distribution and therefore can not be applied directly. One implementable way is to use the supremum likelihood ratio test defined as,
\begin{equation}\label{eq: supQStatLR}%
\underset{\omega\in \Omega}\sup~Q_{L}^{LR}(\omega) = \underset{\omega\in \Omega}\sup ~2\left[ l_n(\hat\beta,\hat\psi, \omega)-l_n(\hat\beta_0,0,\omega)  \right].
\end{equation}

The \textbf{R} package ``\textsf{glarma}" uses the parametrization $\delta^{\prime}=(\beta,\phi,\theta)$ as in \eqref{eq: GLARMA Zt}. For the purpose of fixing the nuisance parameters and optimizing over the remainder the parametrization $\delta=(\beta, \psi,\omega)$ (in \eqref{eq: GLARMA ZtAlt}) is required. Since $\delta= \textrm{A}\delta'$, where $\textrm{A}$ is a fixed nonsingular space matrix of 0's and 1's, the \textsf{glarma} package can be easily modified to optimise the likelihood with respect to $\delta$.

\vskip 8pt
\noindent {\bf 4.2.2. Wald tests}

The Wald test statistic is defined as
\begin{equation}\label{eq: QStatW}%
Q_{L}^{W}=\hat{\psi}^{\mathrm{T}}\hat \Psi^{-1}\hat\psi
\end{equation}
where $\hat\Psi$ is the estimate of the marginal covariance matrix for $\hat\psi$ under the null hypothesis. Under regular cases when there are no nuisance parameters, $\hat\psi = (\hat\phi, \hat\theta)$ and $\hat\Psi$ is the inverse of $I_{L}(\hat \beta_{0})$ in \eqref{eq: InfMat1 GLARMA H0}.

For irregular cases when there are nuisance parameters, the GLARMA\ model
has estimates $\hat\delta=(\hat\beta,\hat\psi,\omega)$, where $\omega$ is fixed.
The supremum Wald test is
\begin{equation}\label{eq: supQStatW}%
\underset{\omega\in \Omega}\sup~Q_{L}^{W}(\omega)= \underset{\omega\in \Omega}\sup ~
\hat{\psi}^{\T} I_{L}(\hat\delta_0)\hat\psi
\end{equation}
where $I_{L}(\hat \delta_{0})$ is defined in \eqref{eq: InfMat GLARMA Red}.

Under the null hypothesis of $\psi=0$, by Theorem \ref{Thm: Aysm GLARMA Nui}, the asymptotic distribution of $Q_{L}^{LR}(\omega)$ and $Q_{L}^{W}(\omega)$ is chi-squared with $L$ (equal to the cardinality of
$J_\phi \bigcup J_\theta$) degrees of freedom for any fixed $\omega$. However, the supremum likelihood ratio and Wald tests are not chi-squared distributed. We will investigate their distributions via simulation.

LRT and Wald tests, under both regular and irregular circumstances, require fitting of the full GLARMA model and, to be justified as a means of screening for serial dependence, they would need to clearly outperform the score test which can be performed using easily applied GLM\ estimates. \par

\section{Simulation}\label{Sec: simulation}%

Simulations are used to assess how well the asymptotic null distribution of supremum score test approximates the finite sample distribution of the supremum score test. We consider the binomial sequences with $m_{t}=2$, $n=200$ and 10,000 replications. The regression is specified with a linear
trend in time and
\begin{equation}\label{eq: Wt simulation}%
W_{t}= -0.5+ (t/n) + Z_{t}
\end{equation}
then the independent sample under the null hypothesis of no serial dependence is simulated by:  $y_{t}|x_{nt}\sim B(m_{t}, 1/(1+\exp(0.5-(t/n))))$. In this section tests are set up againt
GLARMA$(1,1)$ with Pearson residuals as given in \eqref{eq: GLARMA(1,1) ST}.

In Table \ref{tb: QuantileConvergence sup ST} the theoretical quantiles $\mathcal{F}_{\Omega}(u)$ are derived from the distribution in \eqref{eq: sup chi-square1 dist}, which does not require the true
value $\delta_{0}$ thus can be easily implemented in practice. The consistency between theoretical and simulated quantiles starts to break down at the 5\% level on the expanded grid of $[-.99,0.99]$. This can be explained by the requirement for the asymptotic chi-squared distribution of $Q_{1}^{ST}(\omega)$ for fixed $\omega$. The analysis in Section 2.3 implies that $\sup Q_1^{ST}(\omega)$ has asymptotic distribution \eqref{eq: sup chi-square1 dist} only if the $Q_1^{ST}(\omega)$, for fixed $\omega$, is asymptotically chi-squared distributed. Proof of the latter requires $(1-\omega)^{-1}/\sqrt{n}\approx 0$. For $n=200$, $(1-0.99)^{-1}/\sqrt{200}\approx 7.07$. Therefore, for insufficiently large samples, as $\omega\to 1$, the distribution of $Q_1^{ST}(\omega)$ is not chi-squared, which thus affects the distribution for $\sup Q_1^{ST}(\omega)$.

\begin{table}
\caption{Comparison of the consistency of quantiles for approximated theoretical supremum
$\chi^2_{1}$ and simulated supremum test statistics on different scales of the nuisance parameter.}
\label{tb: QuantileConvergence sup ST}\par
\vskip .2cm
\centerline{\tabcolsep=5truept\begin{tabular}{|lccccr|} \hline
&& 10\% & 5\% & 2.5\% & 1\% \\ \hline
\multirow{2}{*}{$\Omega=[-.99,0.99]$}& $\mathcal{F}_{\Omega}(u)$ &5.96& 7.33& 8.69 & 10.51 \\
&$\sup Q_{1}^{ST}$&5.47& 7.73& 11.05& 17.00 \\ \hline
\multirow{2}{*}{$\Omega=[-.80,0.80]$}& $\mathcal{F}_{\Omega}(u)$&4.63 &5.95& 7.29 &9.08\\
&$\sup Q_{1}^{ST}$&4.46& 5.85&7.52& 9.77\\ \hline
\multirow{2}{*}{$\Omega=[-.50,0.50]$}& $\mathcal{F}_{\Omega}(u)$& 3.86 &5.15& 6.45& 8.20\\
&$\sup Q_{1}^{ST}$& 3.81& 5.21 &6.73& 8.58\\ \hline
\end{tabular}}
\end{table}

\section{Applications}\label{Sec: applications}%

\subsection{Example 1: Cambridge-Oxford Boat Race winners -- Bernoulli}

\cite{klingenberg2008regression} modelled the binary time series of the outcome of the Cambridge-Oxford annual boat race with $y_{t}= 1$ when Cambridge wins and $y_{t} = 0$ otherwise. The linear state equation consists of a regression with intercept and the single covariate $x_{t}$ representing the weight difference between the winning and losing side. \cite{klingenberg2008regression} fits a parameter driven model with an AR(1) latent process and his method allows for gaps which occur in the series of 153 race observations over the period 1829 to 2007. Most of the gaps occur early in the series. His fitted model implies the presence of substantial serial dependence and so we use this series as a way of illustrating the performance of the various statistics for detecting serial dependence defined above. Unequal time spacing is not readily accommodated in existing GLARMA modelling software so, in our analysis, time denotes the sequence number of each race.

Simulation results presented here are obtained with 1000 replications of binary sample paths of length $n=153$ generated under the null hypothesis of no serial dependence, $H_{0}: \psi=0$, and using parameter values $\beta_{0}=(0.1937, 0.1176)$, obtained by the GLM fit. All tests are constructed against a GLARMA model with Pearson residuals. Table \ref{tb: Nulldist OxCam Corr} summarizes the upper quantiles of the test statistics against the GLARMA$(1,0)$ model (regular case) and the GLARMA$(1,1)$ model (irregular, nuisance parameter case). The theoretical quantiles for $\mathcal{F}_{\Omega}(u)$ are derived from density \eqref{eq: sup chi-square1 dist}.

The simulated null quantiles of the standard score test $Q_{1}^{ST}(0)$ (in \eqref{eq: STGLARMA Std}) suggest that the limiting chi-squared distribution quantiles slightly overestimate those appropriate for a sample size of $n=153$ in this example.

The quantiles of the supremum score test $\sup Q_{1}^{ST}$ (in \eqref{eq: sup QStatST}) fall below those for the upper tail bound $\mathcal{F}_{\Omega}(u)$ given by \eqref{eq: sup chi-square1 dist} except for the $1\%$ tail probability but again the differences are not substantial, suggesting the upper bound on the tail probabilities of the supremum score statistics provides reasonable guidance on statistical significance in this example.

The simulated quantiles of $\sup Q_{1}^{ST}$, $\sup Q_{1}^{LR}$ from \eqref{eq: supQStatLR} and $\sup Q_{1}^{W}$ from \eqref{eq: supQStatW} are higher than those of $\chi^2_{1}$ as expected, among which, the quantiles of supremum Wald test are substantially higher than supremum score and likelihood ratio tests. A likely explanation for this is that Wald test requires the inverse of the covariance matrix, which is probably poorly estimated with finite samples because of the nuisance parameters.

Results in Table \ref{tb: Nulldist OxCam Corr} indicate that lag 1 serial dependence is significant in the boat race series consistently using the standard score test and the three supremum tests.

\begin{table}[ptb]\centering
\caption{Null distribution upper tail  quantiles and test statistic from the Cambridge-Oxford
boat race series}
\label{tb: Nulldist OxCam Corr}\par
\vskip .2cm
\centerline{\tabcolsep=5truept\begin{tabular}{|llccccr|} \hline
& &10\% & 5\% & 2.5\% & 1\% & Observed \\ \hline
\multirow{2}{*}{against GLARMA$(1,0)$} & $\chi^2_1$ & 2.71 & 3.84 & 5.02 & 6.63 & -\\
& $Q_{1}^{ST}(0)$ & 2.68 & 3.65 & 4.55 & 5.85 &5.69$^\ast$ \\
\multirow{4}{*}{against GLARMA$(1,1)$} & $\mathcal{F}_{\Omega}(u)$ & 5.04 & 6.39 & 7.74 & 9.53 &-  \\
& $\sup Q_{1}^{ST}$ & 4.59& 5.76& 7.72 & 10.82 &11.53$^\ast$ \\
& $\sup Q_{1}^{LR}$ & 5.20& 6.76& 8.57& 11.32 & 10.43$^\ast$ \\
& $\sup Q_{1}^{W}$ & 18.87& 25.01& 32.73& 48.00 &40.65$^\ast$ \\ \hline
\multicolumn{5}{l}{\small $\Omega=-0.9(0.1)0.9$, $\ast$ significant at the $5\%$ level.}
\end{tabular}}
\end{table}

\subsection{Example 2: U.S. Quarterly Recessions -- Binary Response Series}

\cite{kauppi2008predicting} modelled the binary series of U.S. quarterly recessions from 1955:Q4 to 2005:Q4 with a probit link BARMA$(1,1)$ model, and with linear state equation consisting of an intercept and lag 4 interest rate spread, $x_{t-4}$, where lag 4 is selected as a balance between the goodness of in-sample fit and the length of out-sample forecast. Their results indicate that there is clear autocorrelation within the recession series and so it is a good example on which to illustrate the performance of the various test statistics to detect serial dependence, using the logistic BARMA model proposed in Section 3 as the alternative model. When performance on out of sample forecasting with the model is ignored, lag 3 interest rate spread $x_{t-3}$ is found to provide the best within sample fit, and so we use that in the analysis presented here.

To illustrate performance of the likelihood ratio and Wald test BARMA models of appropriate degrees need to be fit to the series. For this series when $L>5$ the estimates of the BARMA models converge slowly or failed to converge and so likelihood ratio and Wald tests are not available for $L>5$. For comparison of these tests with the score and Box-Pearce-Ljung tests we therefore selected $L=5$. Additionally we also used a smaller value $L=3$ to see if serial dependence was effective only at lower lags. Table \ref{tb: ScoreTestsRecessionGLM} shows that all outcomes are significant at the 5\% level and that most serial dependence is contributed by the first 3 lags.

\begin{table}
\caption{Test outcomes for serial dependence in U.S. quarterly recession.}
\label{tb: ScoreTestsRecessionGLM}\par
\vskip .2cm
\centerline{\tabcolsep=5truept\begin{tabular}{|lcr|} \hline
\textbf{Statistic} & \textbf{Value, P-value\thinspace(}$L=3$) & \textbf{Value,
P-value\thinspace(}$L=5$)\\ \hline
Box-Pearce-Ljung: $S_{0}$ & 44.13$(1.42\times10^{-9})$ & 44.25$(2.06\times10^{-8})$\\
BAR\ Score: $S_{1}$ &40.08$(1.02\times10^{-8})$ & 40.59$(1.13\times10^{-7})$\\
BMA\ Score: $S_{2}$ & 36.12$(7.06\times10^{-8})$ & 36.17$(8.78\times10^{-7})$\\
BAR\ LRT: $S_{3}$ & 42.62$(2.97\times10^{-9})$ & 43.57$(2.83\times10^{-8})$ \\
BMA\ LRT: $S_{4}$ & 36.32$(6.40\times10^{-8})$ & 45.66$(1.06\times10^{-8})$ \\
BAR\ Wald: $S_{5}$ & 9.39$(2.45\times 10^{-2})$ & 10.44$(6.38\times 10^{-2})$\\
BMA\ Wald: $S_{6}$ & 18.41$(3.61\times 10^{-4})$ & 39.92$(1.55\times 10^{-7})$\\ \hline
\end{tabular}}
\end{table}%

Next we study the null distributions of the four test statistics. BARMA$(1,2)$ is considered as a pilot example of the alternative model to illustrate the reliability of the asymptotic properties for the four test statistics under finite samples. The study can be generalized to other BARMA models. The simulation is based on 1000 replications of the independent binary series generated with parameter values $\beta_{0}=(-0.223, -1.904)$, obtained by the GLM fit. For score, LRT and Wald statistics, $L=p+q=3$, so they should follow $\chi^{2}_{3}$ distribution asymptotically, and the Box-Pierce-Ljung test $Q_{L}^{BLP}$, $L=3$, is used as a benchmark.

Table \ref{tb: RecessionNulldistquantiles4tests} shows that the quantiles of the Box-Pierce-Ljung have significant upward bias to those of $\chi^2_3$. We believe the most likely explanation for this bias is that the test statistic can, because of the normalization used in residuals, have an extremely large or small value when a binary sequence is dominated by 0's or 1's, which is the situation in this example where there is a high proportion of 0's. Score quantiles are in good agreement with those of $\chi^2_3$ as expected in Theorem \ref{Thm: Asym ST BARMA}, the quantiles of the LRT are upwardly biased similarly to those of the Box-Pearce-Ljung statistic, and, the Wald tests have substantial upward bias. The extremely large quantiles for Wald test are, similarly as in the previous example, due to some of the simulates leading to poorly estimated covariance matrix.
\begin{table}
\caption{Null distribution quantiles of the four test statistics for simulated recession series}
\label{tb: RecessionNulldistquantiles4tests}\par
\vskip .2cm
\centerline{\tabcolsep=5truept\begin{tabular}{|lcccr|}\hline
& 90\% & 95\%  & 97.5\%  &  99\%\\ \hline
$\chi^2_{3}$ & 6.25 & 7.81 & 9.35 & 11.34\\
$Q_{3}^{BLP}$ & 6.95 & 9.74& 11.79& 14.72\\
$Q_{3}^{S}$ & 6.24 & 7.89 & 9.33& 11.18\\
$Q_{3}^{LR}$ &7.54 & 9.63 & 11.31 & 13.99\\
$Q_{3}^{W}$ & 19.96& 32.04 & 53.25 & 71.93\\ \hline
\end{tabular}}
\end{table}

\subsection{Example 3: Court Convictions -- Binomial Response Series}

\cite{dunsmuir2008assessing} considered the number of successful prosecutions obtained from monthly numbers of cases brought to trial in the higher court in the state of NSW, Australia for 6 crime categories: Assault, Sexual Assault, Robbery, Break and Enter, Motor Theft and Other Theft for the period Jan, 1995 to Jun, 2007. Table \ref{tb: CrimeNulldistSuptest} reports summaries of the number of cases brought in each month. There is substantial variation in $m_{t}$ through time in these series and the binomial response distribution for the number of these cases which led to a successful prosecution was used. For each crime, the regressors $X_{t}=(1, T_{t}, \texttt{DNA}_{t-D}, SD_{t})$ are defined as: $T_{t}=t/12$ where $t$ is the month since Jan, 1995; $\texttt{DNA}_{t-D}=\max(t-D-73, 0)$ is a linear increase representing the growth in the number of individuals with DNA records available since Jan, 2001 ($t=73$) when the DNA database was established; and, $D$ is the delay effect of each crime. $SD_{t}$ represents seasonal dummy variables. The logit link was used.

We applied the score tests against the alternative that there is serial dependence of the GLARMA type using Pearson residuals under both regular (no nuisance parameter) and irregular (nuisance parameters present) circumstances, and compare their results. Autocorrelation of Pearson residuals shows that most of the serial dependence is contributed by the first 2 lags and so we used $L=2$ in specifying the score test which is the same form for testing against GLARMA$(0,2)$ or GLARMA$(2,0)$. For each crime, simulation outcomes are obtained with 1000 replications of the independent binomial sequences generated using the GLM\ fit of the real crime data. Table \ref{tb: CrimeNulldistSuptest} gives values for the standard score statistic $Q_{2}^{ST}(0)$ (in \eqref{eq: STGLARMA Std}) to each crime series and when compared with the asymptotic $\chi^2_2$ distribution of Theorem \ref{Thm: QNui GLARMA LimDist} the test statistics is not significant for the crimes of Assault, Motor Theft and Other Theft suggesting that there is no need for serial dependence terms in the linear predictor $W_t$.

We also applied the supremum score test against a GLARMA$(1,2)$ model. In this case there is a one dimensional nuisance parameter. The simulated value of supremum test statistic is obtained as the maximizer of $Q_{2}^{ST}(\omega)$ (in \eqref{eq: QStatST}) over the discrete grid of $\Omega=-0.9(0.1)0.9$ of nuisance parameter values. Table \ref{tb: CrimeNulldistSuptest} shows that the upper tail quantiles of the supremum test are larger than that of chi-squared distribution as is expected. Based on both the standard score test and the supremum score test, the crime categories of Sexual Assault, Break \& Enter and Robbery are serially correlated suggesting strongly that a serially dependent term $Z_t$ should be included in the model for these series.

\begin{table}
\caption{Simulated upper tail quantiles of supremum test statistic applied to crime series of
convictions in the NSW Higher Court}
\label{tb: CrimeNulldistSuptest}\par
\vskip .2cm
\centerline{\tabcolsep=5truept\begin{tabular}{|llllcccccr|} \hline
& \multicolumn{3}{c}{$m_{t}$} & & & & &  & \\
& $\min$ & $\textrm{mean}$ & $\max$ & 10\%& 5\%  & 2.5\% & 1\%  & $Q_{2}^{ST}(0)$
& $\sup Q_{2}^{ST}$ \\ \hline
$\chi^2_2$ & -& -& - & 4.61 &7.38  & 5.99 & 9.21 & - & - \\
\text{Assault} &8 & 85.55& 138 & 5.89 &8.01 &9.39 &13.08 & 2.16 & 7.81\\
\text{SexAssault} & 7& 113.37&250 & 6.15 & 7.72& 8.70 & 12.45 & 39.25$^\ast$ & 42.02$^\ast$\\
\text{BreakEnter} & 4& 53.61 &108 & 5.71 & 7.19 &9.50 &12.98 & 24.17$^\ast$ &42.15$^\ast$\\
\text{Robbery} & 6&90.23 &162 & 5.93 &7.62& 9.29 & 11.06 &24.56$^\ast$ &45.39$^\ast$\\
\text{MotorTheft}&1 & 15.21 &36 & 6.61 & 8.68 & 11.03 & 14.01 &0.06 & 1.06\\
\text{OtherTheft} & 1&10.80 & 31& 6.50 &8.33 &10.20 & 12.47 & 6.98 &7.11\\ \hline
\multicolumn{5}{l}{\small $\Omega=-0.9(0.1)0.9$, $\ast$ significant at $5\%$ level.}
\end{tabular}}
\end{table}

\section{Discussions and Conclusions}\label{Sec: conclusions}%

This paper has developed score tests of the null hypothesis of no serial dependence for time series regression with binomial responses. The test statistics are designed to detect serial dependence of the observation driven type with specific focus on the the GLARMA and BARMA classes of models. Within the GLARMA class, three types of residuals can be specified corresponding to those available in the \texttt{glarma} R-package of \cite{dunsmuir2014glarma}. For BARMA models, except in rather trivial cases when the mean response is constant over time, there are no nuisance parameters occuring under the null hypothesis of no serial dependence and all the test statistics consider (Score, LRT, Wald) are standard with asymptotic chi-squared distribution.

For some GLARMA model specifications nuisance parameters can arise under the null hypothesis and we have demonstrated that the supremum type test of \cite{davies1987hypothesis} can be effective in these situations. Here, we have focussed on a simple situation where there is only one nuisance parameter to deal with but the ideas can be extended to higher dimension for the nuisance parameter space. Implementation of the likelihood ratio, Wald and score tests for the regular case (with no nuisance parameters) can be done easily in the \texttt{glarma} R-package. However to justify the use of these statistics this paper provided the required asymptotic theory primarily via Theorem \ref{Thm: Aysm GLARMA Nui} which extends results of \cite{davis2000ObsDriven} to general regressors but under the null hypothesis of no serial dependence.

For the GLARMA alternatives, the regular score test is of the same form and has the same asymptotic behaviour for both the autoregressive and the moving average specification, that is, it does not discriminate between these two types of alternative dependence.

The simulation results show that the score tests would appear to outperform the Box-Pierce-Ljung, LRT and Wald tests, particularly for binary data which is dominated by 0 or 1. However, more evidence from simulations is needed to investigate the performance of the various statistics under a wider range of models. We have also applied the statistics to various real series to demonstrate the utility of the tests across a diverse range of real settings.

\vskip 14pt
\noindent {\large\bf Supplementary Materials}

The proof of the asymptotic normality of score vectors are standard. Here we present the outline for proof of Theorem \ref{Thm: Aysm GLARMA Nui}. To reduce notation complexity, we give the proof for the example of a GLARMA\ model \eqref{eq: GLARMA ZtAlt} with Pearson residuals. Extension to other type of residuals is straightforward.

For any fixed $\omega$, the true parameter is $\delta_{0} = (\beta_{0},0,\omega)$, and the general parameter, given $\omega$, is $\delta = (\beta, \psi, \omega)$. Let $u = \sqrt{n}(\delta - \delta_{0})$, for these choice the state equation is
\begin{equation*}
W_{t}(\delta)= x_{t}^{T}\beta + Z_{t}(\delta); \quad Z_{t}(\delta) = \underset{j\in J_{\phi}}\sum \omega_{j}Z_{t-j}(\delta) + \underset{J_{\phi}\bigcup J_{\theta}}\sum \psi_{j}e_{t-j}(\delta).
\end{equation*}
Define a linearized version of state equation as
\begin{equation*}
W_{t}^{\dag}(\delta) = x_{t}^{T}\beta + Z_{t}^{\dag}(\delta); \quad Z_{t}^{\dag}(\delta) = \underset{j\in J_{\phi}}\sum\omega_{j}Z_{t-j}^{\dag}(\delta) + \underset{J_{\phi}\bigcup J_{\theta}}\sum\psi_{j}e_{0,t-j}.
\end{equation*}
Following the approach of \cite{davis2000poisson}, linearization is applied to approximate the likelihood function by a convex function of the parameters. In terms of the linearized state equation, let
\begin{equation*}
l_{n}^{\dag}(\delta) = \sum_{t=1}^n \left[ y_t W_{t}^{\dag}(\delta)- m_{t}b(W_{t}^{\dag}) + c(y_{t})\right]
\end{equation*}
where $\delta= \delta_{0} + n^{-1/2}u$, and
\begin{eqnarray*}
R_{n}^{\dag}(u)= - l_{n}^{\dag}(\delta_{0}+ n^{-1/2} u)+ l_{n}^{\dag}(\delta_{0});&
R_{n}(u)= - l_{n}(\delta_{0}+ n^{-1/2} u) + l_{n}(\delta_{0}),
\end{eqnarray*}
It is easy to show that $R_{n}^{\dagger}(u)$ is convex in $u$, the rest of the proof is given in
two major steps:
\begin{enumerate}
\item Establish the limit for $R_{n}^{\dagger}(u)$ as a quadratic form in $u$
plus a normal random variable linear combination of $u$.

\item Show that $R_{n}(u)-R_{n}^{\dagger}(u)\rightarrow0 $ in probability, uniformly for $\left\Vert u\right\Vert <K $ for any finite $K$.
\end{enumerate}

Note $R_{n}^{\dagger}(u)$ can be rewritten as the sum of two parts: $ R_{n}^{\dagger}(u)= B_{n}^{\dag}(u)-A_{n}^{\dag}(u)$, where
\begin{equation*}
A_{n}^{\dag}(u) = \frac{1}{\sqrt{n}} \sum_{t=1}^{n}(y_{t}-m_t\pi_{0,t})%
\begin{bmatrix}
x_{t}\\
\left(1-\underset{j\in J_{\phi}\bigcap J_{\theta}}\sum \omega_{j}\xi^j\right) ^{-1}e_{0,t-J_{\phi}\bigcup J_{\theta}}
\end{bmatrix}^{\T} u = \frac{1}{\sqrt{n}}\sum_{t=1} ^{n}\left( y_{t} - m_t\pi_{0,t}\right)H_{0,t}^{T}u
\end{equation*}
For any fixed $u$, it can be shown using the central limit theorem in \cite{scott1973central} that $A_{n}^{\dag}(u)\overset{d}\rightarrow u^{T} N(0, I(\delta_{0}))$.

Again, for any fixed $u$, there is $\|u^{\ast}\|<\|u\|$ such that
\[
B_{n}^\dag(u)=\frac{1}{2}u^{T}\left(\frac{1}{n}\sum_{t=1}^{n}\sigma_{0,t}^2H_{0,t}H_{0,t}^{T} \right) u + E_{n}^{\dag}(u^{\ast})
\]
as $n\to \infty$, by Chebyshev's inequality
\begin{equation*}
u^{T}\left(\frac{1}{n}\sum_{t=1}^{n}\sigma_{0,t}^2 H_{0,t} H_{0,t}^{T}\right) u \overset{p}\to u^{\T} I_{n}(\delta_0)u
\end{equation*}
where $\underset{n\to\infty}\lim I_{n}(\delta_0)=I(\delta_0)$. Under Condition \ref{cond: Xt},
\begin{equation*}
E_{n}^{\dag}(u^{\ast}) = \frac{1}{6}\|u^{\ast,T}\| ^3 \sum_{t=1}^{n}m_t b^{(3)}(W_{t}(\eta_0 + n^{-1/2} u^\ast))\Vert n^{-1/2} H_{0,t}\Vert ^3 \to 0
\end{equation*}
Then $B_{n}^\dag(u) - u^{T}I(\delta_0)u/2 \overset{p}\rightarrow 0$.

Applying a standard result for functional limit theorems, the $\hat{u}_{n}^{\dag}$ that minimizes $R_{n}^\dag(u)$ satisfies $\hat{u}_{n}^\dag\overset{d}\to \hat{u}^\dag$, where
$\hat{u}^\dag \sim N(0,I(\delta_0)^{-1})$ (see \cite{pollard1991asymptotics}).

Finally we will show that $R_n(u)-R_n^\dag(u)\to 0$. Note by Taylor expansion,
\[
R_{n}(u)-R_{n}^{\dag}(u) = \frac{1}{2}u^{T}\left(\frac{1}{n}\sum_{t=1}^{n}(y_{t}
- m_{t}\pi_{0,t})\ddot{W}_{0,t} \right) u + E_{n}(u^{\ast}) - E_{n}^{\dag}(u^{\ast})
\]
For any fixed $u$, it can be shown using Chebyshev's inequality that the first component converges to zero in probability. $E_{n}^{\dag}(u^\ast)\to 0$ as shown above. There is $\|u^\ast\|< \|u\|$ such that $\delta^\ast = \delta_{0}+u n^{-1/2}$. Variables evaluated at $\delta^\ast$ are, also denoted for example, $\pi_{t}^\ast$,
\begin{equation*}
E_{n}(u^{\ast})=\frac{1}{6}\|u\|^3 l^{(3)}_{n}(\delta_0 + n^{-1/2}u^\ast) =
\frac{1}{6}\|u\|^3n^{-3/2}\sum_{t=1}^{n}\left( (y_{t}-m_{t}\pi_{t}^{\ast})W^{(3),\ast}_{t} - \sigma_{t}^{2,\ast}(4-2\pi_{t}^\ast)\Vert \dot{W}_t^{\ast}\Vert^3 \right)
\end{equation*}
where components of the matrix $W^{(3)}_{t}$ are of the general form
\[
\kappa_{t}=(1-\underset{j\in J_{\phi}\bigcap J_{\theta}}\sum \omega_{j}\xi^{j})^{-1}\frac{\partial^2 e_{t-a}(\delta)}{\partial \delta\partial\delta^{\T}}, \quad a \in J_{\phi}\bigcup J_{\theta}.
\]
By Condition \ref{cond: Xt} and \ref{cond: omega}, for any given $\delta$, as $n\to\infty$,
$\sup_{t\le n} n^{-1/2}\|\dot{W}_{t}\|\to 0$ and $\sup_{t\le n} n^{-1/2}\|\kappa_{t}\|\to 0$.
It follows that $E_{n}(u^\ast)\overset{p}\to 0$.

It is also plausible to conclude that $R_{n}(u)-R_{n}^{\dag}(u)\overset{p}\rightarrow 0$ for
$\Vert u\Vert \le K$, $K<\infty$. Then $\hat{u}_{n}=\arg\min R_n(u)\overset{p}\to \arg\min R_{n}^\dag(u)\overset{d}\to \hat{u}^\dag=\arg\min R^\dag(u)$. And as shown above,
$\hat{u}^\dag\sim N(0,I^{-1}_{\delta_{0}}(\omega))$.

\bibhang=1.7pc
\bibsep=2pt
\fontsize{9}{14pt plus.8pt minus .6pt}\selectfont
\renewcommand\bibname
\end{document}